\documentclass{article}
\usepackage{amsmath}
\usepackage{amsfonts}
\usepackage{amsthm}
\usepackage{eurosym}
\newcommand{\ds}{\displaystyle}

\begin{document}
\title{Distances between the winning numbers in Lottery}
\author{Konstantinos Drakakis}
\date{16 March 2005}
\maketitle

\abstract{We prove an interesting fact about Lottery: the winning 6 numbers (out of 49) in the game of the Lottery contain two consecutive numbers with a surprisingly high probability (almost 50\%).}

\section{Introduction}

The game of lottery exists and has been run in many countries (such as the UK, the US, Germany, France, Ireland, Australia, Greece, Spain, etc.) for a number of years. In this game, the player chooses $m$ numbers from among the numbers $1,\ldots,n>m$, the order of the choice being unimportant and the values of $n$ and $m$ varying from country to country; the lottery organizers choose publicly $m$ numbers in the same way, and if they are the same with the ones the player chose, the player wins. Newspapers usually publish the winning set of numbers along with statistics on the number of times each particular number from $1$ to $n$ has appeared in the winning set. It is however a slightly different and more elusive statistical observation that will be of interest to us here.

Some people have noticed that, in the usual case $m=6$ and $n=49$, it happens very often that at least two of the winning numbers are ``close'' to each other. As $6$ out of $49$ is not really many, this seems at first to be paradoxical, if not altogether wrong, and may remind us strongly of another very similar famous paradox, the Birthday Paradox. In this work we will prove that this observation is well founded, even if we adopt the strictest interpretation of numbers being ``close'', i.e. that they be consecutive. Our problem to solve then will be the following:

{\centering\emph{``What is the probability that, out of $m>0$ numbers drawn uniformly randomly from the range $1,\ldots, n>m$, at least two are consecutive?''}} 

We will calculate this probability in two ways below: one quite ``mechanical'', by finding a recursion and then solving it by means of generating functions, and one combinatorial, which will actually yield a more general result. We will also see that this problem, at least for the usual values $m=6$ and $n=49$, leads to a novel and unexpected gambling application. 

\section{First solution}

Let $f(n,m)$ be the number of ways in which $m$ numbers can be chosen out of $1,\ldots,n$ so that no two are consecutive. For any particular choice, one of the following will hold:
\begin{itemize}
	\item Neither $1$ nor $n$ is chosen: we have to choose $m$ numbers among $2,\ldots,n-1$ and the number of ways this can be accomplished in is $f(n-2,m)$.
	\item $1$ and/or $n$ is chosen: the number of ways this can be accomplished in is, according to the inclusion-exclusion principle, the sum of the number of ways of choosing $1$ and choosing $n$ minus the number of ways in choosing both. Observe now that $2$ cannot be chosen if $1$ is, and that $n-1$ cannot be chosen if $n$ is. Then, in the first two cases the number of choices is $f(n-2,m-1)$, and in the last one $f(n-4,m-2)$, so that the total number of choices if $1$ and/or $n$ is chosen is $2f(n-2,m-1)+f(n-4,m-2)$.
\end{itemize}

Accordingly, summing both cases:
\[f(n,m)=f(n-2,m)+2f(n-2,m-1)-f(n-4,m-2)\]

In addition to the recursive formula above, we need some boundary conditions as well, corresponding to $n=0,1,2,3$ and $m=0,1$. They are provided by the following:
\begin{itemize}
	\item We can choose no numbers in only one way: $f(n,0)=1,\ \forall n\geq 0$.
	\item We can choose one number in $n$ ways: $f(n,1)=n,\ \forall n\geq 0$.
	\item $f(3,2)=1$
\end{itemize}

Let us now write down the generating function for $f(n,m)$:
\[F(z,w)=\sum_{n=4}^\infty \sum_{m=2}^{\left\lceil \frac{n}{2}\right\rceil}f(n,m)z^n w^m\]
The upper boundary for $m$ is determined by the fact that $f(n,m)=0$ if $\ds m\geq \left\lceil \frac{n}{2}\right\rceil+1$. 

By multiplying the recursion formula by $z^nw^m$, and applying the operator $\ds \sum_{n=4}^\infty \sum_{m=2}^{\left\lceil \frac{n}{2}\right\rceil}$, we get:
\[F(n,m)=F_1(n,m)+2F_2(n,m)-F_3(n,m)\]
where 
\begin{eqnarray*}
F_1(n,m)&=&\sum_{n=4}^\infty \sum_{m=2}^{\left\lceil \frac{n}{2}\right\rceil}f(n-2,m)z^n w^m\\
F_2(n,m)&=&\sum_{n=4}^\infty \sum_{m=2}^{\left\lceil \frac{n}{2}\right\rceil}f(n-2,m-1)z^n w^m\\
F_3(n,m)&=&\sum_{n=4}^\infty \sum_{m=2}^{\left\lceil \frac{n}{2}\right\rceil}f(n-4,m-2)z^n w^m
\end{eqnarray*}

For each of the three functions, we get 
\begin{multline*}
F_1(n,m)=\sum_{n=2}^\infty \sum_{m=2}^{\left\lceil \frac{n}{2}\right\rceil+1}f(n,m)z^{n+2} w^m=
z^2 \sum_{n=2}^\infty \sum_{m=2}^{\left\lceil \frac{n}{2}\right\rceil}f(n,m)z^n w^m=
z^2 \left[\sum_{n=4}^\infty \sum_{m=2}^{\left\lceil \frac{n}{2}\right\rceil}f(n,m)z^n w^m+f(3,2)z^3w^2\right]=\\
=z^2 \left[F(z,w)+z^3w^2\right]
\end{multline*}

\begin{multline*}
F_2(n,m)=\sum_{n=2}^\infty \sum_{m=2}^{\left\lceil \frac{n}{2}\right\rceil+1}f(n,m-1)z^{n+2} w^m=
\sum_{n=2}^\infty \sum_{m=1}^{\left\lceil \frac{n}{2}\right\rceil}f(n,m)z^{n+2} w^{m+1}=
w z^2 \sum_{n=2}^\infty \sum_{m=1}^{\left\lceil \frac{n}{2}\right\rceil}f(n,m)z^n w^m=\\
=w z^2 \left[\sum_{n=4}^\infty \sum_{m=2}^{\left\lceil \frac{n}{2}\right\rceil}f(n,m)z^n w^m+f(3,2)z^3w^2+\sum_{n=2}^\infty f(n,1)z^n w\right]=
z^2 w \left[F(z,w)+z^3w^2+w \sum_{n=2}^\infty nz^n\right]
\end{multline*} 

\begin{multline*}
F_3(n,m)=\sum_{n=0}^\infty \sum_{m=2}^{\left\lceil \frac{n}{2}\right\rceil+1}f(n,m-2)z^{n+4} w^m=
\sum_{n=0}^\infty \sum_{m=0}^{\left\lceil \frac{n}{2}\right\rceil}f(n,m)z^{n+4} w^{m+2}=
w^2 z^4 \sum_{n=0}^\infty \sum_{m=0}^{\left\lceil \frac{n}{2}\right\rceil}f(n,m)z^n w^m=\\
=w^2 z^4 \left[\sum_{n=4}^\infty \sum_{m=2}^{\left\lceil \frac{n}{2}\right\rceil}f(n,m)z^n w^m+ f(3,2)z^3w^2+\sum_{n=0}^\infty f(n,1)z^n w+ \sum_{n=0}^\infty f(n,0)z^n\right]=\\
=z^4 w^2 \left[F(z,w)+z^3w^2+w \sum_{n=0}^\infty nz^n+ \sum_{n=0}^\infty z^n\right]
\end{multline*}  

We still need three auxiliary computations:
\begin{itemize}
	\item $\ds \sum_{n=2}^\infty nz^n=z\sum_{n=2}^\infty nz^{n-1}=z\left(\frac{z^2}{1-z}\right)'=z^2\frac{2-z}{(1-z)^2}$
	\item $\ds \sum_{n=0}^\infty z^n=\frac{1}{1-z}$
	\item $\ds \sum_{n=0}^\infty nz^n=z\sum_{n=0}^\infty nz^{n-1}=z\left(\frac{1}{1-z}\right)'=\frac{z}{(1-z)^2}$
\end{itemize}

Putting all of the above together, and after some further algebraic simplifications, we find:
\[F(z,w)= w^2z^4\frac{3+z(z-3+w(z-1)^2)}{(z-1)^2(1-z-wz^2)}\]
Of course, this is not the full generating function, as the cases $n=1,2,3$ and $m=0,1$ are entirely missing; we omitted them in order to avoid to have to deal with ``weird'' boundary conditions such as $f(-3,-1)$ etc. But now we can add them back. 

Remember that $f(n,0)=1,\ n\geq 0$ and $f(n,1)=n,\ n\geq 1$; but we have already carried out the relevant computations as auxiliary computations above. Therefore:
\[\mathcal{F}(z,w)=F(z,w)+z^3w^2+\frac{1}{1-z}+\frac{zw}{(1-z)^2}\]
where the first fraction is the generating function for $f(n,0)$ and the second for $f(n,1)$.
After some algebraic simplifications, we find:
\begin{multline*}
\mathcal{F}(z,w)=\frac{1+zw}{1-z-wz^2}=\frac{1+zw}{1-z-wz^2}=\frac{1}{z}\frac{z(1+zw)}{1-z(1+wz)}=\frac{1}{z}\sum_{n=1}^\infty[z(1+wz)]^n=\\= \sum_{n=1}^\infty\sum_{m=0}^n\left(\begin{array}{c}n\\m\end{array}\right)z^{n+m-1}w^m= \sum_{m=0}^\infty\sum_{n=m}^\infty\left(\begin{array}{c}n-m+1\\m\end{array}\right)z^nw^m
\end{multline*}
so that 
\[f(n,m)=\left(\begin{array}{c}n-m+1\\m\end{array}\right)\]

If then we draw $m$ numbers from the range $1,\ldots, n$, the probability no two are consecutive is:
\[q(n,m)=\frac{\left(\begin{array}{c}n-m+1\\m\end{array}\right)}{\left(\begin{array}{c}n\\m\end{array}\right)}\]
so that the solution to our original problem is:
\[p(n,m)=1-\frac{\left(\begin{array}{c}n-m+1\\m\end{array}\right)}{\left(\begin{array}{c}n\\m\end{array}\right)}\]

We should note here that a proof of the formula for $f(n,m)$ based on induction appears in \cite{R}. 

\section{Second solution}

The second solution, combinatorial in nature, allows us to solve a more general problem: in how many ways $f_k(n,m)$ can we choose $m$ numbers among the numbers $1,\ldots,n$ so that the minimum distance between any two of our choices (which we will be calling the distance of our choice) is $k>0$? There is a very simple formula for that.

Imagine we have numbered $n$ balls with the numbers $1,\ldots,n$, and that we have chosen the numbers $1\leq N_1<\ldots<N_m\leq n$. For every number chosen but the last one, remove the numbers of the $m-1$ balls immediately following it; as for the remaining balls, renumber them consecutively and in the order they are. We will end up with $n-(k-1)(m-1)$ balls numbered consecutively from $1$ to $n-(k-1)(m-1)$, and $(k-1)(m-1)$ blank ones. This final situation will not depend on the balls we chose originally, although the exact positioning of the blank balls among the numbered ones will. Notice finally that the original number of every ball can be recovered: it is the number of balls preceding it, including itself!

Any valid choice of $m$ numbers in the original numbering will correspond to a choice of $m$ numbers after renumbering, and vice versa: after we choose $m$ numbers between $1$ and $n-(k-1)(m-1)$, we insert blanks as described above and renumber, getting a valid choice of numbers in the original numbering. This correspondence is obviously bijective. Therefore,

\[f_k(n,m)=\binom{n-(k-1)(m-1)}{m},\ n>m>1,k\geq 1\]   

For $k=2$ we recover the result of our first solution, and hence the same probability $p(n,m)$ of at least two choices being consecutive. We also obtain the more general formula

\[\ds p_k(n,m)=1-\frac{\left(\begin{array}{c}n-(k-1)(m-1)\\m\end{array}\right)}{\left(\begin{array}{c}n\\m\end{array}\right)}\]

for the probability that at least two of the winning numbers have a distance less than $k$.  

\section{Application in gambling}

The probability $p(n,m)$ can actually be quite large, maybe unexpectedly large: for example, for the usual values $n=49$ and $m=6$, we find $p(49,6)\approx 0.495198$. Therefore, the observation that the winning six numbers of the lottery often contain two that are very close is well founded; in almost one game out of two the winning set of numbers contains two consecutive ones! 

Moreover, as $p(49,6)$ is very close to $0.5$, the problem we just studied can be turned into a successful casino game:  the player bets \euro{e} that $6$ numbers randomly chosen among $1,\ldots,49$ will contain at least two consecutive ones. If this happens, the player gets \euro{e} from the house, otherwise the house wins the player's money. This game is almost fair, as the player has an almost $50\%$ chance to win; but he actually has slightly less than that, and this gives the house a (profitable) advantage!

\section{A slight variant}

What would happen, though, if the player suggests that numbers $1$ and $n$ be treated as consecutive as well, namely if we order the numbers on a ring instead of a line? There should now be fewer possible choices for non-consecutive numbers. Indeed, let now $g_k(n,m)$ be the number of possible choices of $m\leq n$ among $n>0$ numbers so that the minimum distance between any two of the chosen ones is $k$; in other words, among any two chosen numbers, with the property that no number between them is chosen, there are at least $k-1$ numbers lying between them. Then, we can split the choices into those in which one number among $1,\ldots,k-1$ is chosen, and those in which this is not the case: 
\begin{itemize}
	\item If one ball among $1,\ldots,k-1$ is chosen, then the remaining $m-1$ balls can be chosen among $n-2k+1$ balls (we exclude the chosen ball and the $k-1$ adjacent balls on either side); but now, by removing a block of $2k-1$ balls from the circle, we turn it into a line, so the total number of choices, for a fixed choice within $1,\ldots,k$, is $f_k(n-2k+1,m-1)$; and since every different choice within $1,\ldots,k$ leads to different possible choices, the total number of choices in this category is $(k-1) f_k(n-2k+1,m-1)$. 
	\item If no ball is chosen among $1,\ldots,k-1$, then we can just remove them, turn the circle into a line, and renumber: we need to choose $m$ balls among the remaining $n-k+1$, obeying the distance restrictions, and this can happen in $f_k(n-k+1,m)$ ways. Therefore,  
\end{itemize}

\[g_k(n,m)=(k-1) f_k(n-2k+1,m-1)+f_k(n-k+1,m),\ n>m>1,k\geq 0\]	

If we define now
\[\mathfrak{p}_k(n,m)=1-\frac{g_k(n,m)}{\left(\begin{array}{c}n\\m\end{array}\right)}=
1-\frac{\left(\begin{array}{c}n-k+1-(k-1)(m-1)\\m\end{array}\right)+\left(\begin{array}{c}n-2k+1-(k-1)(m-2)\\m-1\end{array}\right)}{\left(\begin{array}{c}n\\m\end{array}\right)}\] 
we find that $\mathfrak{p}_2(49,6)=\mathfrak{p}(49,6)\approx 0.503203$. Therefore, if some casino agreed to play this variant of the game with a player, the player would have a slight advantage over the house, and the latter would loose money!

Table \ref{prob} gives the values of $p_k(49,6)$ and $\mathfrak{p}_k(49,6)$ for $k\in\mathbb{N}^*$:

\begin{table}
\[
\begin{array}{|c||c|c|}
\hline
k& p_k(49,6) & \mathfrak{p}_k(49,6)\\
\hline
\hline
1& 0 & 0\\ \hline
2& 0.495198 & 0.503203\\ \hline
3& 0.766686 & 0.806793\\ \hline
4& 0.903824 & 0.937157\\ \hline
5& 0.966031 & 0.984296\\ \hline
6& 0.990375 & 0.997447\\ \hline
7& 0.99806 & 0.999821\\ \hline
8& 0.999785 & 0.999999\\ \hline
9& 0.999994 & 1\\ \hline
\geq 10 & 1 & 1 \\ \hline
\end{array}
\]
\caption{\label{prob} The probabilities that the winning set of numbers of the standard Lottery has a minimum distance $k$.}
\end{table}

\section{Acknowledgements} 

The author would like to thank an anonymous student of his for communicating to the author his observation about the frequency of appearance of consecutive numbers in the set of the Lottery winning numbers, and thus stimulating him to write this article.

\end{document}